\documentstyle[12pt]{article}

\begin{document}

\begin{titlepage}

\title{Higher analytic torsion of sphere bundles and continuous cohomology of
$\Diff(S^{2n-1})$}

\author{Ulrich Bunke\thanks{Mathematisches Institut, Universit\"at G\"ottingen,
Bunsenstr. 3-5, 37073 G\"ottingen, Germany, bunke@uni-math.gwdg.de}}
\end{titlepage}

\newcommand{\rT}{{\rm T}}
\newcommand{\ch}{{\bf ch}}
\newcommand{\diag}{{\rm diag}}
\newcommand{\proof}{{\it Proof.$\:\:\:\:$}}
\newcommand{\Bbb}{\rm}
\newcommand{\dist}{{\rm dist}}
\newcommand{\kaaa}{{\bf k}}
\newcommand{\paaa}{{\bf p}}
\newcommand{\taaa}{{\bf t}}
\newcommand{\haaa}{{\bf h}}
\newcommand{\R}{{\bf R}}
\newcommand{\Q}{{\bf Q}}
\newcommand{\Z}{{\bf Z}}
\newcommand{\C}{{\bf C}}
\newcommand{\K}{{\tt K}}
\newcommand{\Naaa}{{\bf N}}
\newcommand{\gaaa}{{\bf g}}
\newcommand{\maaa}{{\bf m}}
\newcommand{\aaaa}{{\bf a}}
\newcommand{\naaa}{{\bf n}}
\newcommand{\brr}{{\bf r}}
\newcommand{\res}{{\rm res}}
\newcommand{\Tr}{{\rm Tr}}
\newcommand{\cT}{{\cal T}}
\newcommand{\dom}{{\rm dom}}
\newcommand{\db}{{\bar{\partial}}}
\newcommand{\g}{{\gaaa}}
\newcommand{\cZ}{{\cal Z}}
\newcommand{\cH}{{\cal H}}
\newcommand{\cM}{{\cal M}}
\newcommand{\interi}{{\rm int}}
\newcommand{\singsupp}{{\rm singsupp}}
\newcommand{\cE}{{\cal E}}
\newcommand{\cV}{{\cal V}}
\newcommand{\cI}{{\cal I}}
\newcommand{\cC}{{\cal C}}
\newcommand{\mod}{{\rm mod}}
\newcommand{\cK}{{\cal K}}
\newcommand{\cA}{{\cal A}}
\newcommand{\cEp}{{{\cal E}^\prime}}
\newcommand{\cU}{{\cal U}}
\newcommand{\Hom}{{\mbox{\rm Hom}}}
\newcommand{\vol}{{\rm vol}}
\newcommand{\cO}{{\cal O}}
\newcommand{\End}{{\mbox{\rm End}}}
\newcommand{\Ext}{{\mbox{\rm Ext}}}
\newcommand{\rk}{{\mbox{\rm rank}}}
\newcommand{\im}{{\mbox{\rm im}}}
\newcommand{\sign}{{\rm sign}}
\newcommand{\spann}{{\mbox{\rm span}}}
\newcommand{\symm}{{\mbox{\rm symm}}}
\newcommand{\cF}{{\cal F}}
\newcommand{\Ree}{{\rm Re }}
\newcommand{\Res}{{\mbox{\rm Res}}}
\newcommand{\Imm}{{\rm Im}}
\newcommand{\inter}{{\rm int}}
\newcommand{\clo}{{\rm clo}}
\newcommand{\tg}{{\rm tg}}
\newcommand{\ee}{{\rm e}}
\newcommand{\Li}{{\rm Li}}
\newcommand{\cN}{{\cal N}}
 \newcommand{\conv}{{\rm conv}}
\newcommand{\op}{{\mbox{\rm Op}}}
\newcommand{\tr}{{\mbox{\rm tr}}}
\newcommand{\cs}{{c_\sigma}}
\newcommand{\ctg}{{\rm ctg}}
\newcommand{\degg}{{\mbox{\rm deg}}}
\newcommand{\Ad}{{\mbox{\rm Ad}}}
\newcommand{\ad}{{\mbox{\rm ad}}}
\newcommand{\codim}{{\mbox{\rm codim}}}
\newcommand{\Gr}{{\mbox{\rm Gr}}}
\newcommand{\coker}{{\rm coker}}
\newcommand{\id}{{\mbox{\rm id}}}
\newcommand{\ord}{{\mbox{\rm ord}}}
\newcommand{\nat}{{\bf  N}}
\newcommand{\supp}{{\mbox{\rm supp}}}
\newcommand{\spec}{{\mbox{\rm spec}}}
\newcommand{\Ann}{{\mbox{\rm Ann}}}
\newcommand{\aca}{{\aaaa_\C^\ast}}
\newcommand{\acag}{{\aaaa_{\C,good}^\ast}}
\newcommand{\acage}{{\aaaa_{\C,good}^{\ast,extended}}}
\newcommand{\ck}{{\cal K}}
\newcommand{\tck}{{\tilde{\ck}}}
\newcommand{\tnk}{{\tilde{\ck}_0}}
\newcommand{\ceep}{{{\cal E}(E)^\prime}}
 \newcommand{\ncE}{{{}^\naaa\cE}}
 \newcommand{\Or}{{\rm Or }}
\newcommand{\Diff}{{\cal D}iff}
\newcommand{\cB}{{\cal B}}
\newcommand{\hc}{{{\cal HC}(\gaaa,K)}}
\newcommand{\hcma}{{{\cal HC}(\maaa_P\oplus\aaaa_P,K_P)}}

\newcommand{\vsl}{{V_{\sigma_\lambda}}}
\newcommand{\czg}{{\cZ(\gaaa)}}
\newcommand{\csl}{{\chi_{\sigma,\lambda}}}
\newcommand{\cR}{{\cal R}}
\def\hB{\hspace*{\fill}$\Box$ \newline\noindent}
\newcommand{\varho}{\varrho}
\newcommand{\ind}{{\rm index}}
\newcommand{\Ind}{{\rm Ind}}
\newcommand{\Fin}{{\mbox{\rm Fin}}}
\newcommand{\cS}{{\cal S}}
\newcommand{\orig}{{\cal O}}
\def\hB{\hspace*{\fill}$\Box$ \\[0.5cm]\noindent}
\newcommand{\cL}{{\cal L}}
 \newcommand{\cG}{{\cal G}}
\newcommand{\npci}{{\naaa_P\hspace{-1.5mm}-\hspace{-2mm}\mbox{\rm coinv}}}
\newcommand{\pki}{{(\paaa,K_P)\hspace{-1.5mm}-\hspace{-2mm}\mbox{\rm inv}}}
\newcommand{\mki}{{(\maaa_P\oplus \aaaa_P, K_P)\hspace{-1.5mm}-\hspace{-2mm}\mbox{\rm inv}}}
\newcommand{\Met}{{\cal M}et}
\newcommand{\npi}{{\naaa_P\hspace{-1.5mm}-\hspace{-2mm}\mbox{\rm inv}}}
\newcommand{\ngp}{{N_\Gamma(\pi)}}
\newcommand{\gbg}{{\Gamma\backslash G}}
\newcommand{\gkm}{{ Mod(\gaaa,K) }}
\newcommand{\ggkm}{{  (\gaaa,K) }}
\newcommand{\pkm}{{ Mod(\paaa,K_P)}}
\newcommand{\ppkm}{{  (\paaa,K_P)}}
\newcommand{\makm}{{Mod(\maaa_P\oplus\aaaa_P,K_P)}}
\newcommand{\mmakm}{{ (\maaa_P\oplus\aaaa_P,K_P)}}
\newcommand{\cP}{{\cal P}}
\newcommand{\gm}{{Mod(G)}}
\newcommand{\gk}{{\Gamma_K}}
\newcommand{\La}{{\cal L}}
\newcommand{\cug}{{\cU(\gaaa)}}
\newcommand{\cuk}{{\cU(\kaaa)}}
\newcommand{\dc}{{C^{-\infty}_c(G) }}
\newcommand{\gdk}{{\gaaa/\kaaa}}
\newcommand{\dgkm}{{ D^+(\gaaa,K)-\mbox{\rm mod}}}
\newcommand{\dgm}{{D^+G-\mbox{\rm mod}}}
\newcommand{\vect}{{\C-\mbox{\rm vect}}}
 \newcommand{\cig}{{C^{ \infty}(G)_{K} }}
\newcommand{\gami}{{\Gamma\hspace{-1.5mm}-\hspace{-2mm}\mbox{\rm inv}}}
\newcommand{\cQ}{{\cal Q}}
\newcommand{\mmap}{{Mod(M_PA_P)}}

\newtheorem{prop}{Proposition}[section]
\newtheorem{lem}[prop]{Lemma}
\newtheorem{ddd}[prop]{Definition}
\newtheorem{theorem}[prop]{Theorem}
\newtheorem{kor}[prop]{Corollary}
\newtheorem{ass}[prop]{Assumption}
\newtheorem{con}[prop]{Conjecture}
\newtheorem{prob}[prop]{Problem}
\newtheorem{fact}[prop]{Fact}

\maketitle

\begin{abstract}
Using the higher analytic torsion form of Bismut and Lott
we construct a characteristic class $\cT(E)\in H^*(B,\R)$ for smooth fibre bundles $E\rightarrow B$ with fibre isomorphic to $S^{2n-1}$.
When $E$ is the unit sphere bundle of a hermitean vector bundle 
$V\rightarrow B$, then we calculate this class in terms of the Chern character $\ch(V)$. 

The higher analytic torsion form of the trivial sphere bundle 
over the space of all unit-volume metrics on $S^{2n-1}$
is closed and invariant under the group of diffeomorphisms
$\Diff(S^{2n-1})$. It leads to a continuous cohomology class
$\cT\in H_c^*(\Diff(S^{2n-1}),\R)$. 
Viewing $S^{2n-1}$ as the boundary at infinity of the complex
hyperbolic space $\C H^n$ we obtain an embedding $SU(n,1)\hookrightarrow \Diff(S^{2n-1})$. If $\Gamma\subset SU(n,1)$ is a torsion-free
cocompact subgroup, then the compact K\"ahler manifold 
$B:=\Gamma\backslash \C H^n$ is a model of the classifying space $B\Gamma$. 
Under the identification $H^*(\Gamma,\R)\cong H^*(B,\R)$
the restriction of $\cT$ to $\Gamma$ corresponds to the
characteristic class $\cT(SB)$
of the sphere bundle of $B$.
We show that $0\not= \cT_{4k}\in H^{4k}_c(\Diff(S^{2n-1}),\R)$ as long as $0<2k\le n$.
\end{abstract}


\section{Introduction}
Let  $E\rightarrow B$ be a smooth fibre bundle with a closed manifold  $M$ as fibre, and
let  $F\rightarrow E$ be a complex vector bundle with a flat connection.
In \cite{bismutlott95} Bismut and Lott defined the higher analytic torsion form $T(E,F,T^HE,g^v,h^F)\in\Omega^{ev}(B)$
which depends on the following additional structure:
\begin{itemize}
\item a vertical Riemannian metric $g^v$,
\item a horizontal distribution $T^HE$,
\item a hermitean metric $h^F$ on $F$.
\end{itemize}
Under certain conditions the higher analytic torsion form is closed and defines a cohomology class
$\cT(E,F)\in H^*(B,\R)$ which only depends on $E\rightarrow B$ and $F$ and not on the additional structure.
The class $\cT(E,F)$ behaves natural with respect to pull-back of bundles.

The basic challenge in the theory is to give a topological description of the cohomology class $\cT(E,F)$.
The definition of higher Reidemeister torsion by Igusa and Klein \cite{igusaklein93}
and by Dwyer, Weiss and Williams \cite{dwyerweisswilliams95} yields topologically defined cohomology classes of $B$ with similar  properties as $\cT(E,F)$. The relation between the topologically and analytically
defined classes is far from being understood.

On the computational side the calculation of higher Reidemeister torsion is still very complicated.
In fact explicit results are only known when $E\rightarrow B$ is a quotient
of the Hopf fibration $S^3\rightarrow S^2$ by a finite subgroup $\Z/r\Z$ of $U(1)$ and
$F$ is the bundle $S^3\times_{\Z/r\Z,\rho} \C$, where $\rho:\Z/r\Z\rightarrow U(1)$ is a nontrivial representation.

The higher analytic torsion $\cT(E,F)$ can be calculated in a much larger class of examples.
The structure group $G$ of the fibre bundle $E\rightarrow M$ is a subgroup of the group of diffeomorphisms
$\Diff(M)$ of  $M$. If $G$ is compact, then it was shown by Bismut and Lott \cite{bismutlott95}, that $\cT(E,F,T^HE,g^v,h^F)$ can be expressed in terms of characteristic forms of the $G$-principal bundle $P\rightarrow B$ associated to $E$. Here we assume that $T^H E$ is induced from a connection of $P$,
 $g^v$ comes from a $G$-invariant Riemannian metric on $M$, and 
the hermitean metric $h^{F_M}$ on the restriction $F_M$ of $F$ to $M$  is compatible with $G$, too. 
We see that $\cT(E,F,T^HE,g^v,h^F)$ is a closed form which represents a cohomology class $\cT(E,F)\in H^*(B,\R)$.
This class is independent of the choice of the connection on $P$, the $G$-invariant Riemannian metric on $M$,
and the hermitean metric $h^{F_M}$. But in general it is not an invariant of the pair $(E,F)$ alone.   

The $G$-manifold $M$ and $F_M$ give rise to an invariant power series 
$\rT(M,F_M)\in \hat{I}(G)$ on the Lie algebra of $G$. One obtains $\cT(E,F)$ by applying
the Chern-Weil homomorphism to $\rT(M,F_M)$, i.e.
 \begin{equation}\label{eq1}\cT(E,F)=CW_P(\rT(M, F_M))\ .\end{equation}
In \cite{bunke977} the power series $\rT(M,F_M)$ was explicitly calculated in terms
of a $G$-cell decomposition of $M$. 

In order to assure that $\cT(E,F,T^HE,g^v,h^F)$ is closed one usually assumes that 
$M$ is odd-dimensional and that $F_M$ is acyclic, i.e. $H^p(M,\cF_M)=0$ for all $p\ge 1$,
where $\cF_M$ is the locally constant sheaf of parallel sections of $F_M$.
In the present paper we consider sphere bundles, i.e. bundles
with fibre $S^{2n-1}$. For $n\ge 2$ the bundle  $F_M$ is trivial and never acyclic.
If $F$ is the trivial flat hermitean bundle $E\times\C$ (which we omit in our notation from now on),
then by Lemma \ref{kl} the form $\cT(E,T^HE,g^v)$ is still closed provided the volume of the fibres
induced by $g^v$ is constant. 
By Lemma \ref{lem3} the higher analytic torsion form defines a characteristic class $\cT(E)\in H^{ev}(B,\R)$
for any bundle $E\rightarrow B$ with fibre $S^{2n-1}$ which is independent of the additional choices.

Let $G$ be a compact connected Lie group and $\rho:G\rightarrow U(n)$ be a unitary
representation. Then $G$ acts on the unit sphere $S^{2n-1}$.
Let $T\subset G$ be a maximal torus, $t$ be its Lie algebra, 
$t_\C^*$ be the complexified dual of $t$, 
and $\Lambda\subset t_\C^*$ the weight lattice.
Restriction to $t$ identifies $\hat{I}(G)$ with the space
of Weil-invariant power series  $\hat{I}(T)^W$ on $t$. 
We compute $\rT(S^{2n-1})\in   \hat{I}(T)^W$ in terms of the 
weights $\Delta(\rho)\subset \Lambda$ of the representation $\rho$.
If $\alpha\in \Delta(\rho)$, then $m_\alpha$ denotes the multiplicity of the weight $\alpha$.
We assume that $0\not\in\Delta(\rho)$. 

In order to formulate the result in a compact way we introduce the following formal power series
$Q\in \R[x]$:
\begin{equation}\label{eq0}Q(x):=\sum_{j=1}^\infty \frac{(4j+1)!}{2^{4j}((2j)!)^2}\zeta_R(2j+1)\left(\frac{x}{4\pi^2}\right)^{2j}\ ,\end{equation}
where $\zeta_R$ is the Riemann zeta function.

\begin{theorem}\label{theo1}
$$ \rT(S^{2n-1})\sum_{\alpha\in\Delta(\rho)} m_\alpha Q(\alpha)\ .$$
\end{theorem}

Let $B$ be a  K\"ahler manifold, i.e. $B$ is a Riemannian manifold with
holonomy $U(n)$. Let $SB\rightarrow B$ be the unit sphere bundle.
Combining Theorem \ref{theo1} in the case where $\rho$ is the
standard representation of $U(n)$ on $\C^n$ and Equation (\ref{eq1})
we obtain

\begin{theorem}\label{theo2}
$$\cT(SB)=\sum_{j=1}^\infty \frac{(4j+1)!}{2^{4j}(2j)!}\zeta_R(2j+1)\ch_{4j}(TB)\ ,$$
\end{theorem}
where $\ch(TB)=\ch_0(TB)+\ch_2(TB)+\ch_4(TB)+\dots$ denotes the Chern character
of $TB$.

If we take $B=\C P^n$, then for $p\ge 1$ we have 
$\ch_{2p}(TB)=\frac{(n+1)}{p!} H^p$,
where $H:={\bf c}_1(\cO(1))$. In particular, $\ch_{2p}(TB)\not=0$ for 
$p\le n$. Since $\zeta_R(2j+1)\not=0$ for $j\ge 1$ we see that
 $\cT_{4j}(SB)\not=0$ as long as $0<2j\le n$.

Following the ideas of \cite{bunke978} we can consider the trivial fibre bundle
$$\cE=\Met_1(S^{2n-1})\times S^{2n-1}\rightarrow \Met_1(S^{2n-1})\ ,$$ where $\Met_1(S^{2n-1})$
is the Fr\'echet manifold of all volume one Riemannian metrics on $S^{2n-1}$.
The bundle $\cE$ has a natural horizontal distribution $T^H\cE$ and a tautological
vertical Riemannian metric $g^v$. The group $\Diff(S^{2n-1})$ acts
on this bundle and leaves the additional structure invariant. 
The higher analytic torsion form $T(\cE,T^H\cE,g^v)\in \Omega^{ev}(\Met_1(S^{2n-1})) $
is closed and $\Diff(S^{2n-1})$-invariant. Using that $\Met_1(S^{2n-1})$
is contractible we construct a continuous group cohomology cocycle 
$c_T$ on $\Diff(S^{2n-1})$ which represents a continuous group cohomology class
$\cT\in H^{ev}_c(\Diff(S^{2n-1}),\R)$.

The main goal of the present paper is to prove the following theorem.
\begin{theorem}\label{theo3}
The classes $\cT_{4k}\in H_c^{4k}(\Diff(S^{2n-1}),\R)$,  $0<2k\le n$, are nontrivial.
\end{theorem}

Except for $n=1$ the continuous group cohomology of $\Diff(S^{2n-1})$
has not been calculated yet. Thus Theorem \ref{theo3} provides
a lower bound on the dimensions of the continuous  group cohomology of $\Diff(S^{2n-1})$.
If $n=1$, then explicit generators of $H_c^*(\Diff(S^1),\R)$ are known. Since $S^1$ has a non-trivial
fundamental group we can twist with flat bundles $F\rightarrow S^1$ such that
$\cT_2(F)$ is non-trivial.  We refer to \cite{bunke978} for an explicit calculation of this class.

At the end of this introduction let us sketch the main steps of the proof
of Theorem \ref{theo3}.
We view $S^{2n-1}$ as the boundary at infinity of the complex hyperbolic space
$\C H^n$. We thus obtain an embedding of the isometry group
$SU(n,1)$ of $\C H^n$ into $\Diff(S^{2n-1})$. We then consider a discrete torsion-free
cocompact subgroup $\Gamma\subset SU(n,1)$. We want to show that the restriction to $\Gamma$
of  $\cT_{4k}$ is non-trivial for $0<2k\le n$. 

The quotient $B:=\Gamma\backslash \C H^n$
is a model of the classifying space  $B\Gamma$ and we thus have
a natural identification $H^*(\Gamma,\R)\cong H^*(B,\R)$.
It turns out that the restriction of 
$\cT$ to $\Gamma$ corresponds to the class $\cT(SB)\in H^*(B,\R)$ under this identification.

Now $B$ is a closed locally symmetric K\"ahler manifold and we can apply Theorem \ref{theo2}
in order to express $\cT(SB)$ in terms of the Chern character $\ch(TB)$.
The complex projective space $\C P^n$ is the compact dual symmetric space to $\C H^n$.
Since $\ch_{4k}(T\C P^n)\not=0$ for $0\le 2k\le n$ we conclude by the proportionality 
principle that $\ch_{4k}(TB)\not=0$ for $0\le 2k\le n$, too. Since for $k>0$ the class
$\cT_{4k}(SB)$ is a nonzero multiple of $\ch_{4k}(TB)$ it must be nontrivial.
This finishes the sketch of the proof of Theorem \ref{theo3}.

 \section{Higher equivariant torsion of complex representation spheres}

In this section we specialize the results of  \cite{bunke977} to complex representation spheres
and prove Theorems \ref{theo1} and \ref{theo2}.
Let $G$ be a connected compact Lie group.
Higher equivariant torsion of a closed $G$-manifold $M$ was introduced by Lott \cite{lott94} using a normalization which is slightly different from the one  employed in Bismut and Lott \cite{bismutlott95} and in the present paper. But the results of \cite{lott94} and \cite{bunke977} remain true with the present normalization.
The present normalization behaves even better since the heat kernel expressions
are integrable for small times.

We describe the two nessecary changes in \cite{bunke977}. 
First the definition , Def. 3.1, of the higher equivariant torsion has to be replaced by
$$\rT(M)=\phi \int_0^\infty  \Tr_s \frac{N}{2} (1+2D_t^2)\ee^{D_t^2}\frac{dt}{t}\ ,$$
where we employ the notation introduced in \cite{bunke977} and $\phi$
multiplies the component of homogeneity $k$ by $\left(\frac{1}{2\pi\imath}\right)^{-k}$.
The other point is the computation for $G=S^1$ and $M=S^1$.
Here we employ the formula of Bismut and Lott, \cite{bismutlott95}, Prop. 4.13. 
Let $S^1\cong \R/\Z$ act on $M=\R/\frac{1}{r}\Z$, $r\in\nat$.
Let $\alpha\in (s^1)^*_\C$ be the generator of the weight lattice such that $\alpha(\partial_t)=2\pi\imath $. 
Then 
\begin{equation}\label{eq2}
\rT(M)=\sum_{j=1}^\infty \frac{(4j+1)!}{2^{4j}((2j)!)^2} \zeta_R(2j+1)\left(\frac{r\alpha}{4\pi^2}\right)^{2j}\ .
\end{equation}

Let $G$ be a connected compact Lie group and $\rho:G\rightarrow U(n)$ be an unitary representation. The unit sphere $S^{2n-1}\subset\C^n$ becomes a $G$-manifold.
Let $T\subset G$ be a maximal torus and $t$ be its Lie algebra. By restriction we identify the space of
invariant power series on the Lie algebra of $G$ with the space of Weil invariant
power series $\hat{I}(T)^W$ on $t$. We have to calculate the higher equivariant torsion $\rT(S^{2n-1})\in \hat{I}(T)^W$ 
of the $T$-manifold $S^{2n-1}$. 

By Theorem 1.1 \cite{bunke977} we first must compute the contribution of the one-dimensional $T$-orbits 
to the equivariant Euler characteristic $\chi_T(S^{2n-1})$. Then we have to apply
the homomorphism $T_T$ mapping the target group of the equivariant Euler characteristic to
power series on $t$.  Note that $T_T$ is essentially given by (\ref{eq2}).

We assume that $T$ has no fixed point on $S^{2n-1}$. Let $\Delta(\rho)\subset t^*_\C$ be the set of weights of $\rho$. Then $0\not\in \Delta(\rho)$.
A point $\xi\in S^{2n-1}$ generates a one-dimensional $T$-orbit
iff it is a weight vector.  Let $S_\alpha\subset S^{2n-1}$ be the subset of weight vectors
to $\alpha\in\Delta(\rho)$. All points $\xi\in S_\alpha$  
generate the same orbit type $T/T_\alpha$, where $T_\alpha$ is the stabilizer of $\xi$.
The quotient $T\backslash S_\alpha$ is diffeomorphic to the complex projective space
$\C P^{m_\alpha-1}$, where $m_\alpha$ is the multiplicity of the weight $\alpha$.
Note that $\C P^{m_\alpha-1}$ has Euler characteristic $m_\alpha$.
We conclude that
\begin{equation}\label{eq3}
\chi_T(S^{2n-1})=\sum_{\alpha\in \Delta(\rho)} m_{\alpha} [T/T_\alpha] \ .
\end{equation}

Now $T/T_\alpha\cong S^1$ and using Equation (\ref{eq2}) we obtain  
\begin{equation}\label{eq4}
T_T([T/T_\alpha] )=Q(\alpha)\ ,
\end{equation}
where $Q$ is the formal power series
(\ref{eq0}).
Theorem \ref{theo1} immediately follows from Equations (\ref{eq3}) and (\ref{eq4}). \hB

We now specialize to the case that $G=U(n)$ and $\rho=\id$.
Let $T\subset U(n)$ be the subgroup of diagonal matrices
and identify $t$ with $R^n$ such that $\exp(x)=\diag(\ee^{2\pi\imath x_1},\dots,\ee^{2\pi\imath x_n})$,
where $x=(x_1,\dots,x_n)\in\R^n$.
If $\{\alpha_i\}_{i=1}^n$ denotes the base of $t_\C^*$ which is the dual of the standard base of $t\cong \R^n$
multiplied by $2\pi\imath$, then $\Delta(\rho)=\{\alpha_1,\dots,\alpha_n\}$, and each weight occurs
with multiplicity one. By Theorem \ref{theo1} the homogeneous component of degree $2j$
of $\rT(S^{2n-1})$ is 
\begin{equation}\label{eq5}
\rT_{2j}=\frac{(4j+1)!}{2^{4j}((2j)!)^2}\zeta_R(2j+1)\frac{\sum_{l=1}^n
\alpha_l^{2j}}{(4\pi^2)^{2j}}\ ,
\end{equation}
while the components of odd degree vanish. 

Let now $B$ be a $n$-dimensional K\"ahler manifold and $P\rightarrow B$ be the corresponding
$U(n)$-principal bundle associated to the tangent bundle. If $R$ is the curvature of $TB$,
then the Chern character $\ch(TY)\in H^{ev}(B,\R)$ is represented by 
$$\Tr \exp(-\frac{R}{2\pi\imath})\in\Omega^{ev}(B)\ .$$
We can write $$\ch_{2j}(TY)=CW_P\left(\frac{\sum_{l=1}^n
\alpha_l^{2j}}{(2j)!(4\pi ^2)^{2j}}\right)\ .$$
Thus Theorem \ref{theo2} follows from Equations (\ref{eq5}) and (\ref{eq1}).

\section{Continuous cohomology classes of $\Diff(S^{2n-1})$}

In this section we construct nontrivial continuous cohomology classes $\cT_{4j}\in H_c^{4j}(\Diff(S^{2n-1}),\R)$, $0<2j\le n$, which are related to higher analytic torsion. In fact $\Diff(S^{2n-1})$ is a smooth manifold
and we will use differentiable cochains. For more details
on continuous cohomology we refer to Fuks \cite{fuks84}.

Let $\Met_1(S^{2n-1})$ be the manifold
of all unit volume Riemannian metrics on $S^{2n-1}$.
We consider the trivial fibre bundle $\cE:=\Met_1(S^{2n-1})\times S^{2n-1}\rightarrow \Met_1(S^{2n-1})$
with the obvious horizontal distribution $T^H\cE$. The bundle $\cE$ has a 
tautological vertical Riemannian metric $g^v$. Let $\cT(\cE,T^H\cE,g^v)\in \Omega^{ev}(\Met_1(S^{2n-1}))$
be the higher analytic torsion form of Bismut and Lott.
The following Lemma implies that $\cT(\cE,T^H\cE,g^v)$ is closed.

\begin{lem}\label{kl}
Let $E\rightarrow B$ be a smooth fibre bundle with fibre $S^{2n-1}$,
$T^HE$ be a horizontal distribution and $g^v$ be a vertical Riemannian metric.
If the volume of the fibres is a constant function on $B$, then the higher analytic torsion form
$\cT(E,T^HE,g^v)$ is closed.
\end{lem}
\proof
We cannot refer to Bismut and Lott \cite{bismutlott95}, Cor. 3.26, since $H^*(S^{2n-1},\R)$
is non-trivial in dimensions $0$ and $2n-1$.
But we can employ  Thm. 3.23 (loc. cit).
 
Let $H^*\rightarrow B$ be the bundle of fibrewise harmonic forms.
By Hodge theory this bundle can be identified with
the bundle of fibrewise cohomology and thus aquires a natural flat
connection $\nabla^{H^*}$. As a subbundle of a bundle of Hilbert spaces
it also comes equipped with a natural hermitean metric $h^{H^*}$. Thus we  
can consider the characteristic form $f(\nabla^{H^*},h^{H^*})\in \Omega^{odd}(B)$ introduced in Def. 1.7 (loc. cit.), 
where $f(x)=xe^{x^2}$. 
Since $S^{2n-1}$ is odd dimensional the first term of the r.h.s.
of  Thm. 3.23 (3.119), (loc. cit.),  vanishes and we have 
$$d\cT(E,T^HE,g^v)=-f(\nabla^{H^*},h^{H^*})\ .$$

We now argue that $f(\nabla^{H^*},h^{H^*})=0$.
Note that $H^*$ is nontrivial for $*=0$ and $*=2n-1$.
In both cases we have parallel sections
given by the constant function and the fibrewise volume form,
respectively.  The length of these sections at $b\in B$  
w.r.t. $h^{H^*}$ is just the volume of the fibre of  $E$ over $b$.
Thus $H^*$ is trivialized by sections of constant length and 
$f(\nabla^{H^*},h^{H^*})=0$. \hB

The group $\Diff(S^{2n-1})$ acts on $S^{2n-1}$ and $\Met_1(S^{2n-1})$, 
and thus diagonally on $\cE$.
The additional structure $T^H\cE$ and $g^v$ is $\Diff(S^{2n-1})$-invariant.
We conclude that $\cT(\cE,T^H\cE,g^v)$ is invariant, too.

For $k\in\nat$ we construct a continuous group cohomology cocycle
$c_{T,k}$ of  $\Diff(S^{2n-1})$. The constrcution depends on the choice
of a base point $g_0\in\Met_1(S^{2n-1})$.
Let $q:\Met(S^{2n-1})\rightarrow \Met_1(S^{2n-1})$ be the projection
from the manifold of all Riemannian metrics to the manifold of volume one metrics given by $$q(g):=\frac{g}{\vol_g(S^{2n-1})^{\frac{1}{2n-1}}}\ .$$
The projection $q$ is $\Diff(S^{2n-1})$-equivariant.
 
If $f_0,\dots,f_{k}\in \Diff(S^{2n-1})$, then
we define the map $s(f_0,\dots,f_{k}):\Delta^{k}\rightarrow\Met_1(S^{2n-1})$
by 
$$s(f_0,\dots,f_{k}):= q\left(\sum_{j=0}^{k} t_j f_j(g_0)\right)\ ,$$
where $t_j$ are baricentric coordinates of the $k$-dimensional standard simplex $\Delta^{k}$.
Then we define 
$$c_{T,k}(f_0,\dots,f_{k}):= \int_{\Delta^{k}} s(f_0,\dots,f_{k})^* \cT(\cE,T^H\cE,g^v)\ .$$

\begin{lem}
$c_{T,k}$ is a $\Diff(S^{2n-1})$-invariant continuous (differentiable) 
group cocycle of $\Diff(S^{2n-1})$ and the cohomology class $\cT_{k}\in H^{k}_c(\Diff(S^{2n-1}),\R)$
represented by $c_{T,k}$ is independent of the
choice of the base point $g_0$.
\end{lem}

We leave the simple verification to the interested reader (see also \cite{bunke978}, Sec.2).

\begin{ddd}
We define $\cT:=\cT_1+\cT_2+\dots\in H_c^*(\Diff(S^{2n-1}),\R)$.
\end{ddd}

We now consider $S^{2n-1}$ as the boundary at infinity of the complex hyperbolic space
$\C H^n$. The isometry group $SU(n,1)$ of $\C H^n$ acts on $S^{2n-1}$.
We choose any discrete torsion-free cocompact subgroup $\Gamma\subset SU(n,1)$.
It is a well known theorem of A. Borel  that such subgroups exist.
The group $\Gamma$ acts freely on the contractible space $\C H^n$
and the quotient $B:=\Gamma\backslash \C H^n$ is a model for the classifying space 
$B\Gamma$. In particular, there is a natural identification of $H^*(\Gamma,\R)$
with $H^*(B,\R)$.

\begin{prop}\label{prop1}
Under the above identification the restriction of $\cT$ to $\Gamma$
is equal to $\cT(SB)$.
\end{prop}
\proof
We first make the identification $H^*(\Gamma,\R)\cong H^*(B,\R)$ explicit.
Let $T\in\Omega^j(B)$ be a closed form representing a class $[T]\in H^j(B,\R)$. 
We will construct a group cocycle $C_{T,j}$ which represents the corresponding cohomology class in $H^j(\Gamma,\R)$. 

Let $x_0,\dots,x_j$ be points in $\C H^n$, $0\le j\le 2n$.
We want to construct a map $\sigma(x_0,\dots,x_j):\Delta^j\rightarrow \C H^n$
in an $SU(n,1)$-invariant way. Assume that we already have constructed this map
for all subsets of $j-1$ points. Since $\C H^n$ has negative sectional curvature
there is a unique point $m(x_0,\dots,x_j)$, which has equal distance to all points $x_i$
and  such that this distance is minimal. For any pair of points $y,z\in\C H^n$ let $\gamma(y,z,.):[0,1] \rightarrow \C H^n$ be the unique geodesic from $y$ to $z$. 
If $t\in\Delta^j$ we let $u(t)\in\{0,\dots,j\}$ be the smallest index such that $t_{u(t)}\le t_i$
for all $i\in\{0,1,\dots,j\}$. We define 
\begin{eqnarray*}
\lefteqn{\sigma(x_0,\dots,x_j)(t)}\\&&\hspace{-0.5cm}:=\gamma(m(x_0,\dots,x_j),\sigma(x_0,\dots,\widehat{x_{u(t)}},\dots,x_j)(\frac{t_0}{1-t_{u(t)}},\dots,\widehat{\frac{t_{u(t)}}{1-t_{u(t)}}},\dots,\frac{t_j}{1-t_{u(t)}}),(j+1)t_{u(t)} )\ .\end{eqnarray*}

We choose a base point $x\in\C H^n$.
The cocycle $C_{T,j}$ can now be written as
$$C_{T,j}(f_0,\dots,f_j):=\int_{\Delta^j}\sigma(f_0(x),\dots,f_j(x))^*\tilde{T}, \quad f_0,\dots,f_j\in\Gamma\ ,$$
where $\tilde{T}$ denotes the lift of $T$ from $B$ to $\C H^n$.

The subgroup $SU(n)\subset SU(n,1)$ is a maximal compact subgroup.
Let $g_0$ be the volume one metric on $S^{2n-1}$ which is fixed by $SU(n)$.
We identify the orbit $SU(n,1)(g_0)$ with the symmetric space
$SU(n,1)/SU(n)=\C H^n$ and thus obtain an inclusion $\C H^n\hookrightarrow \Met_1(S^{2n-1})$
which is $SU(n,1)$-equivariant.

The restriction of $T(\cE,T^H\cE,g^v)$ to $\C H^n$ is a $SU(n,1)$-invariant closed form $\tilde{T}$. It induces a closed form $T\in\Omega^*(B)$.
We define the cocycles $c_{T,2k}$ and $C_{T,2k}$ of $\Gamma$ using the base point $g_0$.

\begin{lem}\label{lem2}
The group cohomology cocycles $c_{T,k}$ and $C_{T,k}$ of $\Gamma$
represent the same cohomology class.
\end{lem}
\proof
We define a $\Gamma$-invariant $k-1$-cochain $b$ such that $db=C_{T,k}-c_{T,k}$.
Let $f_0,\dots,f_{k-1}\in\Gamma$. We define a map $a(f_0,\dots,f_{k-1}): [0,1]\times \Delta^{k-1}\rightarrow\Met_1(S^{2n-1})$ by
$$a(f_0,\dots,f_{k-1})(r,t):=r \sigma(f_0(g_0),\dots,f_{k-1}(g_0))(t) + (1-r) s(f_0(g),\dots,f_{k-1})(t)\ .$$
Then we put 
$$b(f_0,\dots,f_{k-1}):=\int_{[0,1]\times\Delta^{k-1}} a(f_0,\dots,f_{k-1})^*\cT(\cE,T^H\cE,g^v)\ .$$
We leave it to the interested reader to check that $b$ has the required properties (see also \cite{bunke978}, Sec.2 for a similar argument). \hB

The restriction of $\cE$ to $\C H^n$ is a trivial $SU(n,1)$-equivariant fibre bundle over $\C H^n$
with $S^{2n-1}$ as fibre equipped with a horizontal distribution and a
vertical Riemannian metric. The quotient by $\Gamma$ is a fibre bundle $\cF\rightarrow B$
with fibres $S^{2n-1}$, horizontal distribution $T^H\cF$ and vertical Riemannian metric $g^v$.

Since the higher analytic torsion form is functorial with respect to pull-back we have
$$T=\cT(\cF,T^H\cF,g^v)\ .$$

In order to finish the proof of the proposition it suffices
to prove the equality of cohomology classes $[T]=\cT(SB)$.
In fact, by Lemma \ref{lem2} the class $[T]$ represents the restriction
of $\cT$ to $\Gamma$ under the identification $H^*(\Gamma,\R)\cong H^*(B,\R)$.

We construct a diffeomorphism of fibre bundles $F:SB\rightarrow \cF$.
Since $SB=\Gamma\backslash S\C H^n$, where $S\C H^n$ is the sphere bundle
of $\C H^n$, and $\cF=\Gamma\backslash(\C H^n\times S^{2n-1})$ w.r.t. the diagonal action
it suffices to produce a $\Gamma$-equivariant diffeomorphism $\tilde{F}:S\C H^n\rightarrow\C H^n\times S^{2n-1}$. Let $v\in S\C H^n$. Then $v$ defines a geodesic ray from the base $x\in\C H^n$  of $v$
which hits the sphere at infinity $S^{2n-1}$ in a point $b$. We define $\tilde{F}(v):=(x,b)$.

On $SB$ we have two sets of additional structures.
On the one hand we have the horizontal distribution $T^H SB$ which is induced from the Levi-Civita
connection and the vertical Riemannian metric $g^v$ induced by the locally symmetric K\"ahler
metric of $B$ normalized such that the fibres have unit volume.
On the other hand we have the integrable horizontal distribution $F^{-1} (T^H\cF)$
and the vertical Riemannian metric $F^*g^v$ induced from $\cF$.

In order to prove the  proposition  remains to prove the following 
lemma which assures that the class $\cT(SB)$
is well defined independent of the choice of additional structures.

\begin{lem}\label{lem3}
Let $E\rightarrow B$ be any smooth fibre bundle with fibre $S^{2n-1}$
and $(T^H E,g^v)$ and $((T^HE)^\prime,(g^v)^\prime)$ two sets of additional structures
such that the fibres have volume one with respect to both, $g^v$ and $(g^v)^\prime$.
Then we have the equality of cohomology classes 
$$[\cT(E,T^HE,g^v)]=[\cT(E,(T^HE)^\prime,(g^v)^\prime)]\ .$$
\end{lem}
\proof
We can connect $T^HE$ with $(T^HE)^\prime$ by a smooth path of horizontal distributions
and $g^v$ with $(g^v)^\prime$ by a smooth path of Riemannian metrics
which induce volume one on the fibres. We thus obtain a horizontal distribution
and a vertical Riemannian metric on $[0,1]\times E\rightarrow [0,1]\times B$
which restrict to the given ones on the boundaries. The associated higher analytic torsion form
is a closed form (by Lemma \ref{kl}) on $[0,1]\times B$ which restricts to $\cT(E,T^HE,g^v)$ and $\cT(E,(T^HE)^\prime,(g^v)^\prime)$ at the boundaries. This observation implies the Lemma. \hB

As we already observed in the introduction by Theorem \ref{theo2}
the class $\cT(SB)_{4k}$ is a non-zero multiple of $\ch_{4k}(TB)$ for $k\in\nat$
and $\ch_{4k}(TY)\not=0$ as long as $0\le 2k\le n$.
 By Proposition \ref{prop1} the classes $\cT_{4k}\in H_c^{4k}(\Diff(S^{2n-1}),\R)$ are nontrivial for $0<2k\le n$. This finishes the proof  of Theorem \ref{theo3}.\hB
\underline{Remarks:}
\begin{itemize}
\item The restrictions of the classes $\cT_{4k-2}$ to $\Gamma$ vanish.
\item Let $\Diff(S^{2n-1})^\delta$ be the diffeomorphism group
with discrete topology. The proof of Theorem \ref{theo3}
shows that the restriction of $\cT_{4k}$ to $\Diff(S^{2n-1})^\delta$ remains nontrivial for $0<2k\le n$.
\end{itemize}

\bibliographystyle{plain}

\end{document}